\documentclass[12pt,a4paper]{article}
\usepackage{amsmath, amssymb, latexsym, color}
\long\def\boxit#1#2{\vbox{\hrule\hbox{\vrule\kern#1\vbox{\kern#1\vbox{#2}\kern#1}\kern#1\vrule}\hrule}}
\def\lloc{L^1_{\rm loc}(\mathbb R^+)}
\def\lumi{(L^1+L^{\infty})(\mathbb R^+)}
\def\lu{L^1(\mathbb R^+)}
\def\lp{L^p(\mathbb R^+)}
\def\lq{L^q(\mathbb R^+)}
\def\li{L^{\infty}(\mathbb R^+)}
\def\lui{L^{1,\infty}(\mathbb R^+)}
\def\r{\mathcal R}

\def\s{\mathcal S}
\def\b{\mathcal B}
\def\nux{\nu_{_X}}

\def\Lx{L^1(\nux)}

\def\Lwx{L_w^1(\nux)}
\def\op#1,#2,#3,{\,#1\,\colon \,#2\to#3\,}
\def\opt#1,#2,#3,#4,{\,#1\in#2\longmapsto#3\in #4\,}
\def\abs#1,{\vert {#1}\vert}
\def\nabs#1,{\left\vert#1\right\vert}
\def\norm#1,{\Vert#1\Vert}
\def\nnorm#1,{\Big\Vert#1\Big\Vert}
\def\nunorm#1,{{\Vert#1\Vert}_\nu}
\def\proof{\noindent {\it Proof.} \ }
\def\endproof{$\hfill\Box$\bigbreak}

\renewcommand{\baselinestretch}{1.15}
 \newtheorem{theorem}{Theorem}[section]
\newtheorem{lemma}[theorem]{Lemma}
\newtheorem{proposition}[theorem]{Proposition}
\newtheorem{corollary}[theorem]{Corollary}
\newtheorem{remark}[theorem]{Remark}
\newtheorem{example}[theorem]{Example}

\newcommand{\R}{\mathbb{R}}

\renewcommand{\title}[1]
{\thispagestyle{empty}
\begin{center}
{\Large \bf #1}
\end{center}}

\newcommand{\authors}[1]
{\begin{center}
\renewcommand{\thefootnote}{\fnsymbol{footnote}}
\setcounter{footnote}{3} {\sc #1 }
\end{center}}

\newcommand{\ack}[1]{\footnote{#1}}

\newcommand{\address}[1]
{\vskip 5ex
\renewcommand{\baselinestretch}{1}
\footnotesize \normalsize
 #1 \\
}

\begin{document}

\title{Optimal domain for the Hardy operator}
  
\authors{
Olvido Delgado\ack{Research partially supported by grant BFM2003-06335-C03-01.} and
Javier Soria\ack{Research partially supported by grants MTM2004-02299 and  
2005SGR00556.\\{\sl  Keywords:} Hardy operator, optimal domain, r.i. space, Lorentz spaces, vector
measures.\\{\sl  MSC2000:} 46E30, 46B25.} }

\bigskip

 {\narrower\noindent \textbf{Abstract.} \small{We study the optimal domain for the Hardy operator considered 
with values in a rearrangement invariant space.  In particular, this domain can be 
represented as the space of integrable functions with respect to a vector 
measure defined on a $\delta$-ring. A precise description is given for the case of 
the minimal Lorentz spaces.}\par}

\bigskip

\section{Introduction}

Let $S$ be the Hardy operator defined by
$$
Sf(x)=\frac{1}{x}\int_0^xf(y)\,dy \ , \ \ x\in(0,\infty) \, ,
$$
for any function $f\in\lloc$.  Let $X$ be a \emph{Banach function ideal lattice}
(abbreviated $BFIL$), i.e., $X$ is  a Banach space
 of real  valued measurable functions on $\R^+$, satisfying that if $g\in
X$ and $|f|\le|g|$ a.e., then $f\in X$ and $\Vert f\Vert_X\le \Vert
g\Vert_X$ (see \cite{BS,krein-petunin-semenov} for further information). For such an $X$,  there is a natural space on which $S$ takes values in
$X$, namely,
$$
[S,X]=\{f\colon\R^+\to\R \textnormal{ measurable},\, S|f|\in X\} \,.
$$
The space
$[S,X]$ is a $BFIL$ itself when endowed with the norm $\Vert
f\Vert_{[S,X]}=\Vert S|f|\,\Vert_X$. Obviously, $S\colon [S,X]\to X$
is continuous. Even more, any $BFIL$ $Y$ such that $S\colon Y\to
X$ is well defined (and so $S$ is continuous, since it is a positive
linear operator between Banach lattices \cite[p.\,2]{lindenstrauss-tzafriri}), is
continuously contained in $[S,X]$. That is, $[S,X]$ is the
\emph{optimal domain} for $S$ (considered with values in $X$) within
the class of $BFIL$.

Similar assertions hold for operators $T$ defined by a positive
kernel $K$ (i.e., $Tf(x)=\int_0^\infty f(y)K(x,y)\,dy$) such that
$T|f|=0$ a.e. implies $f=0$ a.e. This general case has been studied
in \cite{curbera-ricker1,curbera-ricker3}, for $K$ defined on $[0,1]\times[0,1]$, where the authors show  
that the optimal domain $[T,X]$ for $T$, is closely related to the
space   $L^1(\nu_{_X})$ of integrable functions with respect to the
vector measure $\nu_{_X}$, defined by $\nu_{_X}(A)=T(\chi_A)$
(assuming $K$ and $X$ satisfy the minimal conditions for $\nu_{_X}$
to be a vector measure with values in $X$). Indeed, under suitable
additional conditions, both spaces coincide and a precise
description of them is given. The case when $K$ is defined on
$\R^+\times\R^+$ has been studied in \cite{delgado3}. Here, the vector measure
$\nu_{_X}$ associated to $T$ is defined on the $\delta$--ring of the
bounded measurable sets of $\R^+$ (there are classical kernel
operators, like the Hilbert transform, for which $\nu_{_X}$ is not
defined for sets of infinite measure). Again, under
suitable conditions, $[T,X]$ coincides with $L^1(\nu_{_X})$. However, the Hardy operator does not
satisfy these  conditions, and we need to find a different argument 
to describe the space $[S,X]$.

In Section~\ref{odri} we will study several general properties of $[S,X]$ in the   case of {\it rearrangement invariant spaces} $X$ (abbreviated r.i.; that is, if
$g\in X$ and $f$ is equimeasurable with $g$, then $f\in X$ and $\norm
f,_X=\norm g,_X$), and show that the domain is never an r.i. space (Theorem~\ref{noesri}). In Section~\ref{virho}, we prove that $[S,X]$ admits a vector valued integral representation,  and in Section~\ref{optlor} we identify this domain for the minimal Lorentz space $\Lambda_{\varphi}$.
\vspace{.3cm}

\section{Optimal domain and r.i.\ spaces}\label{odri}

We start with a particular case where we are able to identify the domain  for $S$. We observe that $\lui$ is a quasi-Banach r.i. space.

\begin{proposition}\label{dolui}
$
[S,\lui]=\lu
$, with equality of norms.
\end{proposition}

\proof
Recall that $\Vert g\Vert_{\lui}=\sup_{t>0}t\lambda_g(t)$, where $\lambda_g(t)=|\{|g|>t\}|$ is the distribution function of $g$ (see \cite{BS}). Let us prove first the following formula for the distribution function of $Sf$:
If $f\in\lloc$, $f\ge 0$, and $\{Sf>s\}$ has finite measure for all $s>0$, then
\begin{align}\label{distsf}
\lambda_{Sf}(t)=\frac1t\int_{\{Sf>t\}}f(x)\,dx. 
\end{align}

\medskip
\noindent
In fact, since $\{Sf>s\}$ is open and has finite measure, then $\{Sf>s\}=\cup_k(a_k,b_k)$, where $0\le a_k<b_k<\infty$ and these intervals are pairwise disjoint. Moreover, if $a_k\not=0$,
$$
\frac1{a_k}\int_0^{a_k}f(x)\,dx=\frac1{b_k}\int_0^{b_k}f(x)\,dx=s,
$$
and hence, for all cases,
$$
\int_{a_k}^{b_k}f(x)\,dx=\int_{0}^{b_k}f(x)\,dx-\int_{0}^{a_k}f(x)\,dx=s(b_k-a_k).
$$
Thus,
\begin{align*}
|\{Sf>s\}|&=\sum_k(b_k-a_k)=\frac1s\sum_k\int_{a_k}^{b_k}f(x)\,dx\\
&=\frac1s\int_{\cup_k(a_k,b_k)}f(x)\,dx=\frac1s\int_{\{Sf>s\}}f(x)\,dx.
\end{align*}

Using (\ref{distsf}) we now have that if $Sf\in\lui$, $f\ge0$, then

\begin{align*}
\Vert Sf\Vert_{\lui}&=\sup_{s>0}s\lambda_{Sf}(s)=\sup_{s>0}\int_{\{Sf>s\}}f(x)\,dx\\
&=\int_{\{Sf>0\}}f(x)\,dx=\Vert f\Vert_{\lu}.
\end{align*}
Conversely, if $0\le f\in L^1(\R^+)$, then
$\lambda_{Sf}(s)<\infty$ for all $s>0$ and so, the equalities above
hold, i.e., $\Vert f\Vert_{\lu}=\Vert Sf\Vert_{\lui}$.
\endproof

\bigskip

We are going to consider the case of the $\lp$ spaces. It is very easy to show that  $[S,\lu]=\{0\}$. For the other indexes we have the following:

\begin{proposition}\label{lpcont}
 $
\lp\varsubsetneq[S,\lp]
$, $1<p\le\infty$.
\end{proposition}

\proof
Hardy's inequality proves  that $\lp\subset[S,\lp]$. Now, fix $\alpha\in(-1,0)$, and define the unbounded  function  $f_{\alpha}(t)=(1-t)^{\alpha}\chi_{(0,1)}(t)$. Observe that $f_{-1/p}\in\lu\setminus\lp$, $1<p<\infty$. An easy calculation gives,
$$
Sf_{-1/p}(t)=\begin{cases}\displaystyle
      \frac{1-(1-t)^{1-1/p}}{(1-1/p)t},&   0<t<1 \\
  \displaystyle  \frac{p}{p-1}\frac1t,  & \text{ } t\ge 1.
\end{cases}
$$
Therefore, we get the counterexample since $Sf_{-1/p}(t)\in \lq$, for all $1<q\le\infty$. Observe that $f^*_{-1/p}\notin[S,\lp]$ and hence $[S,L^p(\R^+)]$ is not r.i.
\endproof

For a $BFIL$ $X$, if we define 
$$
\Gamma_X=\{f\colon\R^+\to\R \textnormal{ measurable},\, Sf^*\in X\}
\,,
$$
with norm $\Vert f\Vert_{\Gamma_X}=\Vert Sf^*\Vert_X$,  then $\Gamma_X$ is the largest r.i.\ $BFIL$ space contained
in $[S,X]$. In fact, if
$f\in\Gamma_X$, then $S|f|\le Sf^*\in X$ and so $f\in[S,X]$, and if
$Y$ is an r.i.\ $BFIL$ contained in $[S,X]$, then for $f\in Y$ we
have that $f^*\in Y$ and so $Sf^*\in X$, that is $f\in\Gamma_X$.

\begin{proposition}\label{lpnori}Given a $BFIL$ $X$, we have the
following:
\begin{itemize}
\item[(a)] If $S\colon X\to X$, then $X\subset[S,X]$.

\item[(b)] If $X$ is r.i., then $\Gamma_X\subset X\cap[S,X]$.

\item[(c)] If $S\colon X\to X$ and $X$ is r.i., then $\Gamma_X=X$.

\item[(d)] If $X$ is an r.i., the following conditions are equivalent:
\begin{itemize}
\item[(d1)] $\Gamma_X\not=\{0\}$.

\item[(d2)] $\chi_{(0,1)}\in\Gamma_X$.

\item[(d3)] $\chi_{(0,1)}(t)+\frac{1}{t}\chi_{(1,\infty)}(t)\in X$.

\item[(d4)] $(L^\infty\cap L^{1,\infty})(\R^+)\subset X$.
\end{itemize}
\end{itemize}
\end{proposition}

\begin{proof}
(a) is obvious. To prove (b), given $f\in\Gamma_X$, since $f^*\le Sf^*\in X$, then $f^*\in X$ and so $f\in X$. (c) follows from (a), (b), and   the fact that $\Gamma_X$ is the
largest r.i.\ contained in $[S,X]$. Finally, observe that for $f=\chi_{(0,1)}$, we have
$Sf(t)=\chi_{(0,1)}(t)+\frac{1}{t}\chi_{(1,\infty)}(t)$, and the equivalences (d1)-(d4) follow easily. For example, if $g\in(L^\infty\cap L^{1,\infty})(\R^+)$, then $g^*(t)\le C\min(1,1/t)=C(\chi_{(0,1)}(t)+\frac{1}{t}\chi_{(1,\infty)}(t))$. Thus, (d3) implies (d4).
\end{proof}

We observe that we only need $X$ to be an r.i. to prove that (d3) implies (d4). Proposition~\ref{lpcont} shows that the embedding in Proposition~\ref{lpnori}-(a) may
be strict. Let us see now an example of an r.i.\ $BFIL$ space for which the
embedding in Proposition~\ref{lpnori}-(b) is also strict (see also Example~\ref{gamla}).

\begin{proposition}
$\Gamma_{ \lumi}\varsubsetneq \lumi\cap[S, \lumi].
$
\end{proposition}

\begin{proof}
Let us see that $S$ is not bounded on $ \lumi$. In fact, if 
$$
g(t)=\frac1{t\log^2(\frac{e^2}t)}\chi_{(0,1)}(t),
$$
then $g$ is a decreasing function in $\lumi$. Now set $f(t)=g(t-1)\chi_{(1,2)}(t).$ Then, $f^*=g$, $Sf\in \lumi$ (observe that since $f\in L^1$ and it is bounded at zero, then $Sf\in L^{\infty}$), and $Sf^*\notin \lumi$:

$$\Vert Sf^*\Vert_{\lumi}=\int_0^1(Sg)^*(t)\,dt=\int_0^1\frac1{t\log(\frac{e^2}t)}\,dt=\infty.
$$
Hence, we have   shown that 
$\Gamma_{ \lumi}\varsubsetneq\lumi\cap[S,\lumi].
$
 \end{proof}

We are going  to show that Proposition~\ref{lpcont} can be extended to any r.i.\ space:

\begin{theorem}\label{noesri}
 If  $X$ is an r.i.\ $BFIL$ Banach space, and $S:X\rightarrow X$, then $X\varsubsetneq[S,X]$. Hence $[S,X]$ is not r.i. (in fact $[S,X]\not\subset  \lumi$).
 \end{theorem}

\proof
Let us prove that we can  find a function in $[S,X]$ which is not in $ \lumi$, and hence not in $X$ either. We start with the following observation: If $f\ge0$,

\begin{equation}
\label{luli}
f\notin  \lumi\iff {\rm for \ every \ } c>0, \ f\chi_{\{f>c\}}\notin\lu.
\end{equation}

It is clear that if for some $c>0$, $f\chi_{\{f>c\}}\in\lu$, then 

$$
f=f\chi_{\{f>c\}} +f\chi_{\{f\le c\}} \in \lumi.
$$
  Conversely, assume $f=g+h$, $h\in\li$. Take $c=2\Vert h\Vert_{\li}>0$. Then,
$$
f\chi_{\{f>c\}}=(g+h)\chi_{\{g+h>2\Vert h\Vert_{\li}\}}\le(g+h)\chi_{\{|g|>\Vert h\Vert_{\li}\}}\le2|g|.
$$
If $g\in L^1(\R^+)$, then $f\chi_{\{f>c\}}\in L^1(\R^+)$.
\medskip

If $X\subset\lu$, we have that $[S,X]\subset[S,\lu]=\{0\}$, and so, by Proposition~\ref{lpnori}-(a), $X=\{0\}$. Hence, $X\nsubseteq \lu$. Thus, we can find a positive and decreasing function $f\in X$ such that if $F(t)=\int_0^tf(x)\,dx$, then $F$ is strictly increasing and not bounded: take $f_1\in X\setminus \lu$, $f_1$ decreasing (and hence $f_1\ge0$). Choose $f_2\in (L^1\cap L^{\infty})(\R^+)$, decreasing and positive everywhere (e.g. $f_2(t)=(1+t^2)^{-1}$). Note that, since $X$ is an r.i.\ $BFIL$,
$(L^1\cap L^\infty)(\R^+)\subset X$ (see \cite[Theorem II.4.1]{krein-petunin-semenov}) and so
$f_2\in X$. Then $f=f_1+f_2$ satisfies the required conditions. Now take $t_1$=1, and by induction, choose $t_{k+1}>t_k$ satisfying that $F(t_{k+1})=2F(t_k)=2^kF(1)$. We are now going to modify $F$ on each interval $(t_k,t_{k+1})$ in such a way that we obtain a new absolutely continuous, positive and increasing function $G$ satisfying that $F(t)\approx G(t)$, and if $g(t)=G'(t)$, a.e.\ $t>0$, then $g\notin  \lumi$. Hence, $g\in[S,X]$ (observe that $S(g)\approx S(f)\in X$), and $g\notin X$.

On the interval $[0,t_1)$, we set $G(t)=F(t)$. Now we observe the following: since 
$$
\int_{t_k}^{t_{k+1}}f(x)\,dx=F(t_{k})\ge F(t_{k-1})=\int_{t_{k-1}}^{t_{k}}f(x)\,dx,
$$
and $f$ is decreasing, then $t_{k+1}-t_k\ge t_k-t_{k-1}\ge t_2-1$. Therefore, the right triangle $T_k$ determined by the vertices $(t_{k+1}-t_2+1,F(t_{k+1})-F(1))$, $(t_{k+1}, F(t_{k+1})-F(1))$, and $(t_{k+1},F(t_{k+1}))$ (which is congruent to the triangle $T_1$: $(1,F(1))$, $(t_2,F(1))$, and $(t_2, F(2))$) is contained in the right triangle $(t_{k},F(t_{k}))$, $(t_{k+1}, F(t_{k}))$, and $(t_{k+1},F(t_{k+1}))$, for each $k\ge 1$ (observe that $T_k$  has side lengths independent of $k$).

 On the interval $[t_k,t_{k+1}-t_2+1]$, we define $G(t)$ to be the line joining the points $(t_{k},F(t_{k}))$ and $(t_{k+1}-t_2+1,F(t_{k+1})-F(1))$. To define $G$ on the interval $(t_{k+1}-t_2+1,t_{k+1})$ we use the following argument: fix a convex function $h$ on $[1,t_2]$, such that $h(1)=F(1)$, $h(t_2)=F(t_2)$, and $h'(t_2^-)=\infty$ (thus, the graph of $h$ is contained in $T_1$). Now, using the congruence between $T_1$ and $T_k$ (call it $A_k$, so that $A_k(T_1)=T_k$) we translate the graph of $h$ to $T_k$, and define $G(t)$, if $t\in(t_{k+1}-t_2+1,t_{k+1})$, by means of the equality 
 $$
 (t,G(t))=A_k(t-t_{k+1}+t_2,h(t-t_{k+1}+t_2))
 $$
(thus, $G(t)=h(t)$ if $t\in(1,t_2$)).  We observe that $G$ is a continuous, increasing function on $[0,\infty)$. Moreover $G(t)\le F(t)$ since, by concavity,  the graph of $F$ is above the line through the points $(t_k,F(t_k))$ and $(t_{k+1},F(t_{k+1}))$, while $G$ is below that line, by construction. On the other hand, if $t\in(t_k,t_{k+1})$ then 
 $$G(t)\ge G(t_k)=F(t_k)=F(t_{k+1})/2\ge F(t)/2,$$
  and we get the other estimate.

Define now $g(t)=G'(t)$, a.e.\ $t>0$. Let us show that $g\notin  \lumi$: Using (\ref{luli}), if we fix $c>0$, and $k\in\mathbb N$, we can find $s\in (1,t_2)$ such that $g(t)>c$, if $t\in (s,t_2)$ (observe that $g(t_2^-)=G'(t_2^-)=h'(t_2^-)=\infty)$. Then, 
$$\qquad \qquad 
\int_{\{x\in(1,t_{k+1}):g(x)>c\}}g(x)\,dx\ge\sum_{j=2}^{k+1}\int_{s-t_2+t_j}^{t_j}g(x)\,dx=k\int_s^{t_2}h'(x)\,dx\mathrel{\mathop{\longrightarrow\infty}\limits_{k\to\infty}}.$$
\endproof

\begin{remark}{\rm
 We observe that without the hypothesis on $X$, Theorem~\ref{noesri} is false. In fact, as we have proved in Proposition~\ref{dolui}, $
[S,\lui]=\lu
$, which is an r.i. space.
}
\end{remark}

\section{Vector integral representation for the Hardy operator}\label{virho}

\vspace{.3cm}

The representation of a linear operator $T$ between function spaces,
as an integration operator with respect to a vector measure $\nu$,
is always interesting since allows to study the properties of $T$
and its domain through the properties of $\nu$ and the space of
integrable functions with respect to $\nu$. However, this
representation may be not possible. In this section, we give
conditions which guarantee that the Hardy operator $S$ has an
integral representation.

 Associated to $S$ we have the finitely additive set function
$$
A \ \longrightarrow \ \nu(A)=S(\chi_A) \, .
$$
Depending on the family of measurable sets $\r$ on which we define
$\nu$, and the space $X$ where we want $\nu$ to take values,
$\nu\colon\r\to X$ may (or may not) be  a vector measure (i.e., well
defined and countably additive). For instance, if $X=L^1(\R^+)$ no
family of measurable sets $\r$ satisfies that $ \nu\colon\r\to X$ is
a vector measure. Consider another example: the set function
$\nu\colon\b(\R^+)\to(L^1+L^\infty)(\R^+)$, where $\b(\R^+)$ is the
$\sigma$--algebra of all Borel subsets of $\R^+$. This set function
is well defined but it is not a vector measure, since taking
$A_j=[j,j+1)$ we have $\Vert \nu(\cup_{j\ge
k}A_j)\Vert_{L^1+L^\infty}=1$, for all $k$. Then, for any r.i.
$BFIL$  $X$, we have that $\nu\colon\b(\R^+)\to X$ is not a vector
measure, since $X$ is continuously contained in
$(L^1+L^\infty)(\R^+)$ (\cite[Theorem II.4.1]{krein-petunin-semenov}).

We now consider the case when $X$ is a Lorentz space. Recall that for an increasing
concave function $\varphi\colon\R^+\to\R^+$, with $\varphi(0)=0$, the
Lorentz space $\Lambda_\varphi$ is defined by
$$
\Lambda_\varphi=\left\{ f\colon\R^+\to\R \textnormal{ measurable}, \
\int_0^\infty f^*(t)d\varphi(t)<\infty\right\} \, ,
$$
where $f^*$ is the decreasing rearrangement of $f$. The space $
\Lambda_\varphi$ endowed with the norm $\Vert f\Vert_{
\Lambda_\varphi}=\int_0^\infty f^*(t)\,d\varphi(t)$, is an r.i.\
$BFIL$ space. Choosing $\r$ as the $\delta$--ring (ring closed under
countable intersections)
\begin{equation}\label{EQ: dring}
\r=\{A\in \b(\R^+):\ |A|<\infty \ \textnormal{ and } \ \exists\,
\varepsilon>0\,,  \ |A\cap[0,\varepsilon]|=0\} \, ,
\end{equation}
where $|\,\cdot |$ is the Lebesgue measure on $\R^+$, we have the
following result.

\begin{proposition}\label{vectmeas}
$\ \nu(A)\in\Lambda_\varphi$ for
every $A\in\r$ if and only if
\begin{equation}\label{EQ: thetaX}
\theta_\varphi(y)=\int_y^\infty\frac{\varphi'(t)}{t}\,dt<\infty \, , \
\textnormal{ for all } y>0\, ,
\end{equation}
where $\varphi'$ is the derivative of $\varphi$. Moreover, if
\eqref{EQ: thetaX} holds, then $\nu\colon\r\to\Lambda_\varphi$ is a
vector measure.
\end{proposition}

\noindent {\it Proof.} We first observe that (\ref{EQ: thetaX}) is equivalent to saying that $\theta_{\varphi}$ is integrable near 0, since 
$$
\int_0^{\varepsilon}\theta_{\varphi}(y)\,dy=\varphi(\varepsilon)-\varphi(0^+)+\varepsilon\theta_{\varphi}(\varepsilon).
$$
Now, given $A\in\r$ we have
$$
\int_0^\infty
\nu(A)^*(t)\,d\varphi(t)=\varphi(0^+)\,\nu(A)^*(0^+)
+\int_0^\infty \nu(A)^*(t)\,\varphi'(t)\,dt \, ,
$$
where
$$
\nu(A)^*(0^+)=\norm \nu(A),_\infty=
\sup_{0<x<\infty}\,\frac{1}{x}\int_0^x \chi_A(y)\, dy =
\sup_{0<x<\infty}\frac{1}{x}\,|[0,x]\cap A|\le1,
$$
and since $(S\abs f,)^*\le Sf^*$,
\begin{eqnarray*}
\int_0^\infty \nu(A)^*(t)\,\varphi'(t)\,dt & \le &
\int_0^\infty \frac{\varphi'(t)}{t}\int_0^t\chi_{[0,|A|)}(y)\,
dy\, dt \\ & = & \int_0^{|A|}
\int_y^\infty\frac{\varphi'(t)}{t}\, dt\, dy \, .
\end{eqnarray*}
Then, if  \eqref{EQ: thetaX} holds, 
$\nu(A)\in\Lambda_{\varphi}$, for all $A\in\r$. 

Conversely, if $\nu(A)\in\Lambda_\varphi$ for
every $A\in\r$, then, taking $A=[\frac{a}{2},a]$ for any $a>0$ we
have $A\in\r$ and
$$
\frac{a}{2}\,\theta_\varphi(a)\le\int_{\frac{a}{2}}^a\theta_\varphi(y)dy=
\int_0^\infty\nu(A)(t)\varphi'(t)dt\le
\int_0^\infty\nu(A)^*(t)\varphi'(t)dt<\infty \, ,
$$
since $\theta_\varphi$ is decreasing. So, $\theta_\varphi(y)<\infty$
for all $y>0$. Hence, $\varphi$
satisfying \eqref{EQ: thetaX} is equivalent to
$\nu\colon\r\to\Lambda_\varphi$ is well defined. Let us see that in
this case $\nu$ is countably additive:

Given a disjoint
sequence $(A_j)\subset\r$, with $A=\cup_{j\ge1}A_j\in\r$, and taking
$\varepsilon>0$ such that $|A\cap[0,\varepsilon]|=0$, we have
$$
\sup_{0<x<\infty}\,\frac{1}{x}\,|[0,x]\cap \cup_{j\ge
k}A_j|\le\frac{1}{\varepsilon}\,|\cup_{j\ge k}A_j| \, .
$$
Then
$$
\norm \nu(\cup_{j\ge k}A_j),_{\Lambda_{\varphi}} \le
\frac{\varphi(0^+)}{\varepsilon}\,|\cup_{j\ge k}A_j|
+\int_0^{|\cup_{j\ge k}A_j|} \theta_{\varphi}(y)\, dy \
\longrightarrow \ 0
$$
as $k\to\infty$, since $|A|<\infty$ and condition (\ref{EQ: thetaX}) holds. $\hfill \Box$

\bigskip

From Proposition~\ref{vectmeas} we deduce conditions for a general space $X$,
under which $\nu\colon\r\to X$ is a vector measure. Let $X$ be an
r.i.\ $BFIL$ space and $\varphi_{_X}$ the fundamental function of $X$
defined by $\varphi_{_X}(t)=\Vert \chi_{[0,t]}\Vert_X$, for
$t\in\R^+$. Taking an equivalent norm in $X$ if necessary, we have
that $\varphi_{_X}$ is concave (\cite{BS, krein-petunin-semenov}). Then, since
$\Lambda_{\varphi_{_X}}$ is continuously contained in $X$ (see \cite[Theorem II.5.5]{krein-petunin-semenov}), we have that a measure with values in
$\Lambda_{\varphi_{_X}}$ is also a measure with values in $X$.

\begin{corollary}
If $\varphi_{_X}$ satisfies
\eqref{EQ: thetaX}, then $\nu\colon\r\to X$ is a vector measure.
\end{corollary}

\vspace{-0.05cm}

\begin{remark}{\rm
 If $X$ has fundamental function $\varphi_{_X}$ satisfying \eqref{EQ:
thetaX} and $\varphi_{_X}(0^+)=0$, it is sufficient to take
$\tilde{\r}=\{A\in\b(\R^+):\, |A|<\infty\}$ for $\op
\nu,\tilde{\r},X,$ to be a vector measure. }
\end{remark}

From now on we will assume that $X$ is an r.i.\
$BFIL$, with fundamental function $\varphi_{_X}$ satisfying
\eqref{EQ: thetaX}.  Thus, $\op
\nu,\r,X,$ is a vector measure, which will be denoted by $\nux$ to
indicate the space where the values are taken. We will make use of the integration theory for vector measures defined on $\delta$-rings,
due to Lewis \cite{lewis2} and Masani and Niemi
\cite{masani-niemi1, masani-niemi2}. So, we consider the
space $\Lx$ of integrable functions with respect to $\nux$, namely,
measurable functions $\op f,\R^+,\R,$ such that
\begin{itemize}\setlength{\leftskip}{-2ex}\setlength{\itemsep}{.5ex} 
\item[(i)] $f$ is integrable with
respect to $\abs x^*\nux,$, for all $x^*\in X^*$, and

\item[(ii)] for each $A\in\b(\R^+)$, there is a vector, denoted by $\int_A
fd\nu\in X$, such that
\begin{equation*}
x^*\left(\int_Afd\nu\right)=\int_A fdx^*\nu, \ \hbox{ for all }
x^*\in X^*,\ 
\end{equation*}
\vspace{-.15cm}
\end{itemize}
where $\abs x^*\nux,$ is defined on $\b(\R^+)$ as the variation of
the real measure $x^*\nux$. Noting
that $|A|=0$ if and only if $\nu(A)=0$ a.e., the space $\Lx$ endowed with the norm
$$
\norm f,_{\nux}=\sup_{x^*\in B_{X^*}}\int \abs f, d\abs x^*{\nux},,
$$
is a $BFIL$ space, in which the
$\r$--simple functions (i.e., simple functions with support in $\r$)
are dense. Moreover, $\Lx$ is order continuous (i.e., order bounded
increasing sequences are norm convergent). Since $X$ is a Banach
lattice and $\nux$ is a positive vector measure, it can be proved
that $\norm f,_{\nux}=\norm {\int \abs f,d\nux},_X$, for all
$f\in\Lx$ (see the discussion after the proof of \cite[Theorem
5.2]{curbera-ricker1}). For results concerning the space $L^1$ of a
vector measure defined on a $\delta$--ring, see \cite{delgado2}.

\bigskip

For every $f\in\Lx$ it can be proved that $Sf=\int f d\nux\in X$,
see \cite[Proposition 3.1.(b)]{delgado3}. Thus, $S$ coincides on
$\Lx$ with the integration operator with respect to $\nux$ and
$\Lx\hookrightarrow [S,X]$, with $\norm f,_{[S,X]}=\norm f,_{\nux}$.
Even more,
$L^1(\nu_{_X})$ is the largest order continuous $BFIL$ space contained
in $[S,X]$. Let us prove this fact: Let $Y$ be an order continuous
$BFIL$ such that $Y$ is continuously contained in $[S,X]$. Given $0\le f\in Y$, there are simple
functions $\psi_n$ such that $0\le\psi_n\uparrow f$. We take the
$\r$--simple functions $\varphi_n=\psi_n\chi_{[\frac{1}{n},n]}$ for
which $0\le\varphi_n\uparrow f$. For all $A\in\b(\R^+)$ we have
$0\le\varphi_n\chi_A\uparrow f\chi_A\in Y$. Since $Y$ is order
continuous it follows that $\varphi_n\chi_A\to f\chi_A$ in $Y$ and
then $\varphi_n\chi_A\to f\chi_A$ in $[S,X]$. So $\norm
S(f\chi_A)-S(\varphi_n\chi_A),_X=\norm {S\abs
f\chi_A-\varphi_n\chi_A,\,},_X\to0$ as $n\to\infty$. Thus,
$S(\varphi_n\chi_A)=\int_A\varphi_nd\nux$ converges in $X$, for every
$A\in \b(\R^+)$. Using  \cite[Proposition 2.3]{delgado2}, we have
that $f\in\Lx$. Therefore $Y\subset\Lx$ and the inclusion is positive and
continuous.

If $X$ is order continuous, then it  is easy to see that $[S,X]$ is also
order continuous, and thus $\Lx=[S,X]$.

\bigskip

Now, let us consider the larger space
$$
\Lwx=\left\{\op f,\R^+,\R, \textnormal{ measurable}:\, \int \abs
f,d\abs x^*\nux,<\infty \textnormal{ for all } x^*\in X^*\right\},
$$
which is a $BFIL$ space with the norm
$\Vert\cdot\Vert_{\nu_{_X}}$, satisfying the Fatou property  (i.e., $(f_n)\subset \Lwx$, $\sup_n\norm
f_n,_{\nux}<\infty$, $0\le f_n\uparrow f$ a.e.\ implies $f\in \Lwx$
and $\norm f_n,_{\nux}\uparrow\norm f,_{\nux}$). Note that
$\Lx\hookrightarrow\Lwx$.

In a similar way to \cite[Proposition 3.2.(ii)]{curbera-ricker3}, it
can be proved that $[S,X]\hookrightarrow\Lwx$ with $\norm
f,_{\nux}\le\norm f,_{[S,X]}$. Even more, $L_w^1(\nu_{_X})$ is the smallest
$BFIL$ space with the Fatou property containing $[S,X]$.

If $X$ has the Fatou property, then $[S,X]$ also has the Fatou
property and thus $\Lwx=[S,X]$.

\bigskip

Summarizing, the following result has been established.

\begin{proposition}\label{PROP: all results}
Let $X$ be an r.i.\
$BFIL$ space whose fundamental function $\varphi_{_X}$ satisfies \eqref{EQ:
thetaX}. For the $\delta$--ring $\r$ given in \eqref{EQ: dring} we
have:
\begin{enumerate}\setlength{\leftskip}{-2ex}\setlength{\itemsep}{.8ex}
\item[(a)] $\op \nux,\r,X,$ is a vector measure, where
$\nux(A)=S(\chi_A)$.

\item[(b)] $\Lx\hookrightarrow [S,X] \hookrightarrow \Lwx$.

\item[(c)] $L^1(\nu_{_X})$ is
the largest order continuous $BFIL$ space contained in $[S,X]$.

\item[(d)] $L_w^1(\nu_{_X})$ is the smallest $BFIL$ space with the Fatou
property containing $[S,X]$.

\item[(e)] If $X$ is order continuous, then $\Lx=[S,X]$.

\item[(f)] If $X$ has the Fatou property, then $\Lwx=[S,X]$.
\end{enumerate}
\end{proposition}

\smallskip

\begin{example}{\rm
For $1<p\le\infty$, the space $X=L^p(\R^+)$ satisfies the
hypothesis of Proposition~\ref{PROP: all results}. Since for
$1<p<\infty$ the space $L^p$ is order continuous and has the Fatou
property, we have
$$
[S,L^p]=L^1(\nu_{_{L^p}})=L_w^1(\nu_{_{L^p}}) \ .
$$
For $p=\infty$ we have
$$
L^1(\nu_{_{L^\infty}})\hookrightarrow
[S,L^\infty]=L_w^1(\nu_{_{L^\infty}}) \ ,
$$
since $L^\infty$ has the Fatou property. Observe that
$L^1(\nu_{_{L^\infty}})\varsubsetneq [S,L^\infty]$. For instance,
$\chi_{\R^+}\in [S,L^\infty]\backslash L^1(\nu_{_{L^\infty}})$.
Indeed, if $\chi_{\R^+}\in L^1(\nu_{_{L^\infty}})$, then by
\cite[Corollary 3.2.b)]{delgado2}, $\nu_{_{L^\infty}}$ is strongly
additive (i.e., $\nu_{_{L^\infty}}(A_n)\to0$ whenever $(A_n)$ is a
disjoint sequence in $\r$), but taking $A_n=[2^n,2^{n+1})$ we obtain
$\norm {\nu_{_{L^\infty}}(A_n)},_{\infty}={1}/{2}$, for all
$n\ge1$ and this is a contradiction.}
\end{example}

\smallskip

\begin{example}{\rm 
Let $X$ be a Lorentz space $\Lambda_\varphi$ with
$\varphi$ satisfying \eqref{EQ: thetaX}; that is, satisfying the
hypothesis of  Proposition~\ref{PROP: all results}. Since
$\Lambda_\varphi$ has the Fatou property, we have
$$
L^1(\nu_{_{\Lambda_\varphi}})\hookrightarrow[S,\Lambda_\varphi]=L_w^1(\nu_{_{\Lambda_\varphi}})
\ .
$$
In the case when $\varphi(0^+)=0$ and $\varphi(\infty)=\infty$ we
have that $\Lambda_\varphi$ is order continuous (see \cite[Corollary
1 to Theorem II.5.1]{krein-petunin-semenov}) and so
$$
L^1(\nu_{_{\Lambda_\varphi}})=[S,\Lambda_\varphi]=L_w^1(\nu_{_{\Lambda_\varphi}})
\ .
$$}
\end{example}

\vspace{.5cm}

\section{Optimal domain for the Lorentz spaces $\mathbf \Lambda_{\varphi}$}\label{optlor}

\vspace{.3cm}

Let $X$
be a $BFIL$ space. Recall the definition of  the space
$$
\Gamma_X=\{\op f,\R^+,\R, \textnormal{ measurable},\ Sf^*\in X\} \,
. \vspace{.15cm}
$$
 In
general, $\Gamma_X$ is not a closed subspace of $[S,X]$. For
instance, if we take $X=L^p$ for $1<p<\infty$, we have (see Proposition~\ref{lpcont}):
$$
\s(\r)\subset\Gamma_{L^p}= L^p \varsubsetneq
[S,L^p]=L^1(\nu_{_{L^p}}) \, ,
$$
where $\s(\r)$ is the space of $\r$--simple functions.
Then, $\Gamma_{L^p}$ is not closed in $[S,L^p]$, since $\s(\r)$ is
dense in $L^1(\nu_{_{L^p}})$. 
\bigskip

\begin{example}\label{gamla}{\rm
Consider the Lorentz space $\Lambda_\varphi$. For any measurable
function $f$, noting that $Sf^*$ is decreasing, it follows
\begin{eqnarray*}
\int_0^\infty (Sf^*)^*(t)\,d\varphi(t) & = & \int_0^\infty
Sf^*(t)\,d\varphi(t) \\ & = & \varphi(0^+)Sf^*(0^+)+ \int_0^\infty
Sf^*(t)\,\varphi'(t)\,dt \\ & = & \varphi(0^+)\Vert
Sf^*\Vert_\infty+\int_0^\infty
\frac{\varphi'(t)}{t}\int_0^tf^*(s)\,ds\,dt
\\ & = & \varphi(0^+)\Vert
f\Vert_\infty + \int_0^\infty f^*(s)
\int_s^\infty\frac{\varphi'(t)}{t}\,dt\,ds \\ & = &
\varphi(0^+)\Vert f\Vert_\infty + \int_0^\infty
f^*(s)\,\theta_\varphi(s)\,ds \, .
\end{eqnarray*}
Therefore,
$$
\Gamma_{\Lambda_\varphi}=L^\infty\cap\Lambda_{\int_0^t\theta_\varphi(s)ds}
\, .
$$

In the case when $\varphi(0^+)=0$, we have $
\Gamma_{\Lambda_\varphi}=\Lambda_{\int_0^t\theta_\varphi(s)ds}$.
Moreover, in this case, $ \Gamma_{\Lambda_\varphi}=\Lambda_\varphi$
if and only if  $\int_0^t\theta_\varphi(s)\,ds$ and $\varphi$ are
equivalent (e.g. $\varphi(t)=t^{1/p}$, for $1<p<\infty$), and this
holds if and only if there exists a constant $C>0$ such that
\begin{equation}\label{EQ: phi-constant}
t\,\theta_\varphi(t)\le C\,\varphi(t), \ \ \textnormal{ for all }
t\in(0,\infty) \, ,
\end{equation}
since
\begin{eqnarray*}
\int_0^t\theta_\varphi(s)\,ds & = &
\int_0^t\int_s^\infty\frac{\varphi'(y)}{y}\,dy\,ds= \int_0^\infty
\frac{\varphi'(y)}{y}\int_{[0,t]\cap[0,y]}ds\,dy \\ & = &
\int_0^\infty \frac{\varphi'(y)}{y}\,\min\{t,y\}\,dy =
\int_0^t\varphi'(y)\,dy+t\int_t^\infty\frac{\varphi'(y)}{y}\,dy \\ &
= & \varphi(t)+ t\,\theta_\varphi(t) \ .
\end{eqnarray*}

Condition \eqref{EQ: phi-constant} is also equivalent to saying that $\varphi'\in B_1$ (see \cite{CRS}).

The function $\varphi(t)=\min\{1,t\}\,$ (for which
$\Lambda_\varphi=L^1+L^\infty$) does not satisfy   condition
\eqref{EQ: phi-constant}, so $\Gamma_{L^1+L^\infty}\varsubsetneq
L^1+L^\infty$. (For more information about this kind of embeddings and the boundedness of the Hardy operator see \cite{CRS}.)
}
\end{example}

\bigskip

Now we will   describe the space $[S,\Lambda_\varphi]$ in the case when $\varphi(0^+)=0$. Observe
that
\begin{eqnarray*}
\int_0^\infty (S\abs f,)^*(t)\,\varphi'(t)\,dt & \ge & \int_0^\infty
S\abs f,(t)\,\varphi'(t)\,dt=\int_0^\infty
\frac{\varphi'(t)}{t}\int_0^t\abs f(s),\,ds\,dt \\ & = &
\int_0^\infty \abs f(s), \int_s^\infty\frac{\varphi'(t)}{t}\,dt\,ds
= \int_0^\infty \abs f(s),\,\theta_\varphi(s)\,ds \, .
\end{eqnarray*}
Then, we always have that 
\begin{equation}\label{inclusion}
[S,\Lambda_\varphi]\hookrightarrow
L^1({\theta_\varphi}(t)\,dt), 
\end{equation}
where $L^1({\theta_\varphi}(t)\,dt)$ denotes the
space of integrable functions with respect to the Lebesgue measure
with density $\theta_\varphi$.

\bigskip

We will use the following result for an r.i.\ $BFIL$ $X$, with the Fatou
property. In this case, $X'$ (the K\"{o}the dual
of $X$) is a norming subspace of $X^*$, that is
$$
\Vert f\Vert_X=\sup_{g\in B_{X'}}|<g,f>|=\sup_{g\in
B_{X'}}\Big|\int_0^\infty g(x)f(x)\,dx\Big| \, ,
$$
\cite[Proposition 1.b.18]{lindenstrauss-tzafriri}. Note that if $f$ is positive, the supremum
above can be taken for positive functions in $B_{X'}$. 

\begin{lemma}\label{LEM: weightedL1}
Let $X$ be an r.i.\ $BFIL$ space, with the Fatou property. Suppose $X$ satisfies
\begin{equation}\label{EQ: hy-function}
h_y\in X \ \  { a.e. \ } y>0 \ , \ \  { where } \
h_y(x):=\frac{1}{x}\,\chi_{_{[y,\infty)}}(x) \, .
\end{equation}
Then $\,L^1({\phi_{_X}}(t)\,dt)\hookrightarrow[S,X]\,$, for
$\phi_{_X}(y)=\norm h_y,_X$.
\end{lemma}

\noindent {\it Proof.} \ Note that, since $X$ is and r.i., from
Proposition~\ref{lpnori}-(d) we have that condition \eqref{EQ: hy-function}
is equivalent to $\Gamma_X\not=\{0\}$, and this happens if and only
if $(L^1\cap L^\infty)(\R^+)\subset [S,X]$, since $\Gamma_X$ is the
largest r.i.\ $BFIL$ contained in $[S,X]$. In particular, any simple
function $f$ with finite support is in $[S,X]$ and
\begin{eqnarray*}
\norm f,_{[S,X]} & = & \norm {S\abs f,},_X=\sup_{0\le g\in
B_{X'}}\int_0^\infty g(x)\,S\abs f,(x)\,dx \\ & = & \sup_{0\le g\in
B_{X'}}\int_0^\infty \frac{g(x)}{x}\int_0^x\abs f(y),\,dy\,dx \\ & =
& \sup_{0\le g\in B_{X'}}\int_0^\infty \abs
f(y),\int_y^\infty\frac{g(x)}{x}\,dx\,dy \\ & \le &
\int_0^\infty\abs f(y),\,\norm h_y,_X\,dy=\int_0^\infty\abs
f(y),\,\phi_{_X}(y)\,dy\, .
\end{eqnarray*}
For $f\in L^1(\phi_{_X}(t)\,dt)$ we can take simple functions
$(f_n)$ with finite support, such that $0\le f_n\uparrow |f|$. Then
$$
\sup_{n\ge1} \norm f_n,_{[S,X]}\le \sup_{n\ge1}\int_0^\infty\abs
f_n(y),\,\phi_{_X}(y)\,dy=\int_0^\infty\abs
f(y),\,\phi_{_X}(y)\,dy<\infty \, .
$$
Thus, $f\in[S,X]$ and $\norm f,_{[S,X]}=\sup_{n\ge1}\norm
f_n,_{[S,X]}\le\int_0^\infty\abs f(y),\,\phi_{_X}(y)\,dy$. We have
used that $[S,X]$ has the Fatou property since $X$ has this
property.$\hfill \Box$

\bigskip
\bigbreak

\begin{remark}\label{remb}{\rm (a) If $X$ is an r.i.\ $ BFIL$ space, with
fundamental function satisfying (\ref{EQ: thetaX}), then we have that
$\s(\r)\subset [S,X]$. In particular, $S\chi_A\in X$ for $A=(a,b)$,
with $0<a<b<\infty$. Then, since
$S\chi_A(x)=(1-\frac{a}{x})\chi_{(a,b)}(x)+(b-a)\frac{1}{x}\chi_{[b,\infty)}(x)$
and $(1-\frac{a}{x})\chi_{(a,b)}(x)\in (L^1\cap
L^\infty)(\R^+)\subset X$, condition \eqref{EQ: hy-function} holds
for $X$.

(b) Let $X=\Lambda_\varphi$, with $\varphi$ satisfying (\ref{EQ: thetaX}) and
$\varphi(0^+)=0$. From (a) we have that $h_y\in\Lambda_\varphi$ and
$$
\phi_{_{\Lambda_\varphi}}(y)=\int_0^\infty
h_y^*(s)\,\varphi'(s)\,ds=\int_0^\infty\frac{\varphi'(s)}{y+s}\,ds
\, . \vspace{.2cm}
$$
Actually, in this case, (\ref{EQ: thetaX}) and \eqref{EQ: hy-function} are
equivalent. Then, by Lemma~\ref{LEM: weightedL1},
$\,L^1(\phi_{_{\Lambda_\varphi}}(t)\,dt)\hookrightarrow[S,\Lambda_\varphi]\,$.
Note that $\phi_{_{\Lambda_\varphi}}$ is equivalent to the function
given by $\theta_\varphi(t)+\frac{\varphi(t)}{t}$.\, Indeed,
$$
\phi_{_{\Lambda_\varphi}}(t)=
\int_t^\infty\frac{\varphi'(s)}{t+s}\,ds+\int_0^t\frac{\varphi'(s)}{t+s}\,ds
$$
where
\begin{eqnarray*}
&
\displaystyle\frac{1}{2}\,\theta_\varphi(t)=\frac{1}{2}\int_t^\infty\frac{\varphi'(s)}{s}\,ds\le
\int_t^\infty\frac{\varphi'(s)}{t+s}\,ds\le\int_t^\infty\frac{\varphi'(s)}{s}\,ds=\theta_\varphi(t)
& \\  \\
&
\displaystyle\frac{1}{2}\,\frac{\varphi(t)}{t}=\frac{1}{2t}\int_0^t\varphi'(s)\,ds\le\int_0^t\frac{\varphi'(s)}{t+s}\,ds
\le\frac{1}{t}\int_0^t\varphi'(s)\,ds=\frac{\varphi(t)}{t}. & \\
\end{eqnarray*}
 So, $\phi_{_{\Lambda_\varphi}}(t)\le
\theta_\varphi(t)+\frac{\varphi(t)}{t}\le2\phi_{_{\Lambda_\varphi}}(t)$.
}
\end{remark}
 
\bigskip

\begin{theorem}\label{dominlor}
A Lorentz space $\Lambda_\varphi$ with $\varphi$ satisfying
\eqref{EQ: thetaX}, $\varphi(0^+)=0$ and for which there exists a
constant $C>0$ such that
\begin{equation}\label{EQ: theta-condition}
\frac{\varphi(t)}{t}\le C\,\theta_\varphi(t), \  \ \ { for\ 
all \ \ \ } t\in(0,\infty) \ ,
\end{equation}
satisfies
$$
[S,\Lambda_\varphi]=L^1({\theta_\varphi}(t)\,dt)=L^1({\phi_{_{\Lambda_\varphi}}}(t)\,dt)
\, .
$$
\end{theorem}

\vspace{.2cm}

\proof  Using \eqref{inclusion} and Lemma~\ref{LEM: weightedL1}, we have that 
$L^1({\phi_{_{\Lambda_\varphi}}}(t)\,dt)\hookrightarrow
[S,\Lambda_\varphi]\hookrightarrow L^1({\theta_\varphi}(t)\,dt)$. If
\eqref{EQ: theta-condition} holds, then $\theta_\varphi$ is equivalent to $\theta_\varphi(t)+
\varphi(t)/t$, which is equivalent (by Remark~\ref{remb}-(b)) to $\phi_{_{\Lambda_\varphi}} $.
So,
$
L^1({\theta_\varphi}(t)\,dt)=L^1({\phi_{_{\Lambda_\varphi}}}(t)\,dt)=[S,\Lambda_\varphi]
\, .
$
\endproof

 We consider now the special case of the Lorentz spaces $L^{p,q}$. We show that for $q=1$, the domain coincides with an $L^1$-space with respect to an absolutely continuous measure, but this result does not hold if $1<q\le\infty$:

\begin{proposition} \label{nohaypeso}(a) For $1<p<\infty$,
\begin{equation}\label{iguallp1}
[S,L^{p,1}]=L^1(t^{-1/p'}dt) \, .
\end{equation}

(b) If $1<p<\infty$ and $1\le q\le\infty$, then $L^1(t^{-1/p'}dt)\subset[S,L^{p,q}].$ 

\medskip

(c) For every $1<q\le \infty$, there does not exist a nonnegative function $v\in \lloc$ for which $[S,L^{p,q}]=  L^1(v(t)\,dt).$
\end{proposition}

\proof To prove (a),  we observe that the function $\varphi(t)=t^{1/p}$ satisfies
\eqref{EQ: theta-condition}:
$$\theta_{\varphi}(t)=\frac{1}{p-1}t^{-(1-1/p)}=\frac{1}{p-1}\frac{\varphi(t)}{t}.$$
The result follows from Theorem~\ref{dominlor},  since $\Lambda_{\varphi}=L^{p,1}$

\medskip

(b) is a consequence of (a) and the fact that $L^{p,1}\subset L^{p,q}$. 

\medskip

Suppose now that $[S,L^{p,q}]=  L^1(v(t)\,dt).$ Then, using a small modification of the result in \cite[p.\ 316]{kaak}, it follows that, since $L^1(v(t)\,dt)\subset [S,L^{p,q}]$, there exists a constant $C>0$ such that $C\le t^{1/p'}v(t)$,  and hence $L^1(v(t)\,dt)\subset [S,L^{p,1}]$. Therefore,  $[S,L^{p,q}]=[S,L^{p,1}]$. But, taking a decreasing function $f\in L^{p,q}\setminus L^{p,1}$, we find that $f\in L^{p,q}\subset[S,L^{p,q}]$, and $f\le Sf\in L^{p,1}$, which is a contradiction.
\endproof

\begin{remark}{\rm
Proposition~\ref{nohaypeso} shows that $L^1(t^{-1/p'}dt)$ is the largest $L^1$-space contained in $[S,L^{p,\infty}]$. If we consider the converse embedding $[S,L^{p,\infty}]\subset L^1(v(t)\,dt)$, then a necessary condition is that 
\begin{equation}\label{intevp}
\int_0^{\infty}\frac{v(t)}{t^{1/p}}\,dt<\infty.
\end{equation}
On the other hand, if (\ref{intevp}) holds,
then any decreasing function in $[S,L^{p,\infty}]$ belongs also to $L^1(v(t)\,dt)$.
}
\end{remark}

\address{
\noindent
Olvido Delgado\\ Dept.\  of Mathematics\\ University of Sevilla\\
E-41080 Sevilla, SPAIN \ \ \ {\sl E-mail:} 
 {\tt olvido@us.es}

\medskip

\noindent
Javier Soria\\ Dept.\  Appl.\  Math.\  and Analysis
\\ University of Barcelona\\ E-08007 Barcelona,
 SPAIN\ \ \ {\sl E-mail:} 
 {\tt soria@ub.edu}
}
\end{document}